\documentclass[11pt]{article}
\usepackage{amsmath}
\numberwithin{equation}{section}
\usepackage{amsfonts,amssymb,amsmath,amsthm}
\usepackage{mathrsfs}
\newcommand{\beq}{\begin{equation}}
\newcommand{\enq}{\end{equation}}

\newtheorem{Theorem}{Theorem}[section]
\newtheorem{Lemma}[Theorem]{Lemma}

\newtheorem{Definition}[Theorem]{Definition}
\newtheorem{Remark}[subsection]{Remark}

\newcommand{\benu}{\begin{enumerate}}
\newcommand{\beqa}{\begin{eqnarray}}
\newcommand{\beqan}{\begin{eqnarray*}}
\newcommand{\eay}{\end{array}}
\newcommand{\edm}{\end{displaymath}}
\newcommand{\eenu}{\end{enumerate}}
\newcommand{\eeq}{\end{equation}}
\newcommand{\eeqa}{\end{eqnarray}}
\newcommand{\eeqan}{\end{eqnarray*}}

\newcommand{\br}{\begin{Remark}}
\newcommand{\er}{\end{Remark}}

\newcommand{\bqa}{\begin{eqnarray}}
\newcommand{\eqa}{\end{eqnarray}}
\newcommand{\bqw}{\begin{eqnarray*}}
\newcommand{\eqw}{\end{eqnarray*}}

\newcommand{\bea}{\begin{array}{cc}}
\newcommand{\ena}{\end{array}}

\allowdisplaybreaks[4]
\begin{document}
\begin{center}

{\large \bf Uniform attractor of a non-autonomous Lam\'{e} thermoelastic system}\\

 \vspace{0.20in}Yuming Qin$^{1,2\ast}$ $\ $ Hongli Wang$^{2}$\\
\end{center}
 $^{1}$ Insitute of Nonlinear Science, Donghua University,
Shanghai, 201620,  P. R. China, \\
$^{2}$ School of Mathematics and Statistics, Donghua University,
Shanghai, 201620, P. R. China.
 \vspace{3mm}

\begin{abstract}
	In this paper, we investigate the dynamical behavior of non-autonomous Lam\'{e} thermoelastic systems within 
    $N$-dimensional materials. With appropriate constraints on nonlinear characteristics and functional parameters, we initially establish the existence of a uniformly absorbing set by constructing a Lyapunov function. Subsequently, we employ the contraction mapping principle to demonstrate the uniformly asymptotic compactness of the system. Finally, under irrotational conditions, we prove the existence of a uniform attractor $\mathcal{A}_{\Sigma}$ in the space $H_c$.
 \vskip 3mm


 \smallskip
 \smallskip
 
\end{abstract}

\section{Introduction}
\setcounter{equation}{0}

\let\thefootnote\relax\footnote{*Corresponding author.}
\let\thefootnote\relax\footnote{E-mails: yuming@dhu.edu.cn (Y. Qin), whl@mail.dhu.edu.cn (H. Wang).}
\quad
In this paper, we are concerned with the following non-autonomous Lam\'{e} thermoelastic system
\begin{equation}\label{eq1.1}
	\left\{\begin{array}{l}
		 u_{tt}-\Delta_e u+\alpha(t) \nabla \theta+f(u)=0, \\
		\theta_t-\kappa(t) \Delta \theta+\alpha(t) \operatorname{div} u_t=g(x,t),
	\end{array} \quad \text { in } \quad \Omega \times[\tau,+\infty),\right.
\end{equation}
with initial conditions and boundary conditions
\begin{equation}\label{eq1.2}
	\left\{\begin{array}{l}
		u(x, \tau)=u_\tau(x),\, u_t(x, \tau)=v_\tau(x),\, \theta(x, \tau)=\theta_\tau(x),\quad x\in\Omega,\,\tau\in\mathbb{R}, \\
		u(x,t)=0,\,\theta(x, t)=0,\,x\in\partial\Omega,\,t\geq\tau,
	\end{array}\right.
\end{equation}
where $\Omega \subset \mathbb{R}^N, N \geq 3$, is a bounded domain with smooth boundary $\partial \Omega$ and $\Delta_e$ denotes the Lam\'{e} operator defined by
$$
\Delta_e=\mu \Delta u+(\lambda+\mu) \nabla \operatorname{div} u \text {, with } \mu>0 \text { and } \lambda \in \mathbb{R} \text {. }
$$

In 1974, Achenbach\cite{JA11} first introduced the Lam\'{e}-Navier equation $\rho \vec{f}+(\lambda+G) \nabla(\nabla \cdot \vec{x})+G \nabla \cdot \nabla \vec{x}=\overrightarrow{0}$, which is fundamental in the study of linear elasticity. Since then, the Lam\'{e} system has been recognized by numerous scholars for its significant role in partial differential equations. It serves as the theoretical basis for describing linear elastic media deformation and is central to mathematical physics and elasticity theory. Recently, researchers have focused on exploring the long-time behavior of the Lam\'{e} systems.

For the autonomous Lam\'{e} system, Yamamoto\cite{KY12} studied the autonomous Lam\'{e} system $u_{t t}-\mu \Delta u-(\mu+\lambda) \nabla \operatorname{div} u=0$ under Dirichlet boundary conditions, establishing the exponential decay of solutions through the scattering theory of Lax and Phillips. Then Astaburuaga and Char\~{a}o \cite{AB13} investigated the autonomous Lam\'{e} system $u_{t t}-a^2 \Delta u-\left(b^2-a^2\right) \nabla \operatorname{div} u+q(x) g\left(u_t\right)=0$ within a three-dimensional bounded domain, considering the influence of infinite memories. They obtained a comprehensive and accurate assessment of the solutions' convergence to zero as $t$ approaches infinity. For the non-autonomous Lam\'{e} system, Costa, Freitas and Wang\cite{AL14} researched the long-time dynamics of the non-autonomous Lam\'{e} system $u_{t t}-\mu \Delta u-(\mu+\lambda) \nabla \operatorname{div} u+\alpha u_t+f(u)=b(t)$, to establish the existence and upper-semicontinuity of pullback attractors, which are critical for understanding the system's asymptotic behavior. Subsquently, Bezerra and Nascimento\cite{FD15}  studied a class of semilinear thermoelastic systems with variable thermal coefficients, further advancing the understanding of long-term dynamics in such a system
$$
\left\{\begin{array}{l}
	\partial_t^2 u-\mu \Delta u-(\lambda+\mu) \nabla \operatorname{div} u+\beta_{\varepsilon}(t) \nabla \theta=f(u), \\
	\partial_t \theta-\Delta \theta+\beta_{\varepsilon}(t) \operatorname{div} \partial_t u=0,
\end{array}\right.
$$
they obtained the existence, regularity and upper-semicontinuity of the pullback attractors. It has been observed that these studies focus on the global attractors and pullback attractors of the Lam\'{e} system. Regarding elastic systems, \cite{YQ32} provides a detailed discussion on recent results and asymptotic behaviors of solutions to non-classical thermo(visco)elastic models. Further research in this area can be found in \cite{YQ27,YQ28,YQ29}.

The complexity increases significantly when studying non-autonomous systems. A primary method in this field is the uniform attractor, initially introduced by Haraux \cite{AH17,AH18} and later improved in subsequent studies \cite{VV19,VV20,SV22}. For nonlinear and complex systems, uniform attractors provide a method for understanding and predicting their complex dynamical behavior, which is crucial for comprehending the system's global properties and long-time behavor. Recently, Yang\cite{LY16}  investigated the dynamics of non-autonomous plate-type evolutionary equations $u_{t t}+a(x) u_t+\Delta^2u+\lambda u+f(u)=g(x, t)$ with critical nonlinearity, to further establish the existence of a uniform attractor in the space $H_0^2(\Omega) \times L^2(\Omega)$. In subsequent work, Xie, Li and Zeng\cite{YX21} studied nonclassical diffusion equations with memory $u_t(t)-\Delta u_t(t)-\Delta u(t)-\int_0^{\infty} k(s) \Delta u(t-s) d s  +f(u(t))=g(t)$ and successfully established the existence and structure of a compact uniform attractor. Later, Xiong and Sun\cite{YX24} considered the non-autonomous damped wave equation $\partial_t^2 u+\gamma \partial_t u-\Delta u+f(u)=g(x, t)$, and  derived an upper bound for the $\varepsilon$-entropy of the uniform attractor. 

Building on the aforementioned research, we have fostered a particular interest in the uniform attractors of system (\ref{eq1.1})-(\ref{eq1.2}). We will outline the primary challenges and novel contributions of this study in detail.

(i) It is widely recognized that in the field of thermoelasticity for N-dimensional materials, the total energy may not decay to zero, as demonstrated in \cite{CM25}. Drawing on the works \cite{JB26,FD1} and employing Helmholtz decomposition from the fundamental theorems of vector calculus, we establish the existence of uniform attractor within both the irrotational component of the displacement vector field and the thermal discrepancy.

(ii) A Lyapunov function, as detailed in equation (\ref{eq3.27}), has been constructed to confirm a uniform absorbing set for the system. While previous studies, including \cite{YW23}, have considered a broad range for the constants involved, our work precisely defines the specific bounds for each parameter. These detailed bounds are thoroughly presented in the proof of Lemma \ref{le3.9}.

(iii) The nonlinear term $f(u)$ and the non-autonomous external force term $g(x,t)$ have been considered in system (\ref{eq1.1})-(\ref{eq1.2}), resulting in the thermoelastic system being a highly complex non-autonomous dynamical system. This refined model delivers a more precise representation of the dynamical responses of materials under the concurrent influences of thermal and mechanical forces. However, these two elements also introduce significant challenges in studying the long-time behavior of the solutions.

(iv) The majority of existing studies such as \cite{FD1,FD2,AL14} are dedicated to examining the long-time behavior of solutions for the Lam\'{e} system, including the existence of global attractors, pullback attractors, and their upper semi-continuity and regularity. However, to the best of our knowledge, there are currently no results regarding the uniform attractors for the non-autonomous Lam\'{e} thermoelastic system. Thus, our paper represents a novel attempt.
 
The structure of this paper is as follows. In Section 2, we present the necessary assumptions, and based on \cite{FD1}, drawing upon the research presented in \cite{FD1}, we establish the well-posedness of the system given by equations (\ref{eq1.1})-(\ref{eq1.2}),  as in Theorem \ref{th2.4}. Furthermore, in Section 3, through a series energy estimates, we demonstrate the existence of a uniformly absorbing set in the space $H_c$, specifically for the case of a curl free part, as in Theorem \ref{th3.11}. In Section 4, we derive uniform asymptotic compactness of the corresponding family of processes $\left\{U_G(t, \tau)\right\}(G \in \Sigma)$, generated by (\ref{eq2.7})-(\ref{eq2.8}), as in Theorem \ref{th4.4}. As a result, we obtain the existence of a uniform attractor $\mathcal{A}_{\Sigma}$, as in Theorem \ref{th4.5}.

\section{Assumptions and well-posedness}

 \quad  In this section, we establish the local and global well-posedness of the problem (\ref{eq1.1})-(\ref{eq1.2}) in $H_c$ and $H_d$
 under the assumptions (\ref{eq2.1})-(\ref{eq2.6}).
 
\textbf{\large {2.1 Assumptions}}

As usual, $L^q (\Omega), (1 \geq q \leq \infty)$ and $H^i(\Omega), (i = 1, 2, 3)$, denote the standard
notations of Lebesgue integral and Sobolev spaces. The $L^2$-inner product is denoted by $(\cdot, \cdot)$ and $\|\cdot\|_B$ denotes the norm in the space $B$. For simplicity, we use $\|u\|$ instead
of $\|u\|_2$ when $q = 2$.

We assume that $\alpha, \kappa: \mathbb{R} \rightarrow \mathbb{R}$ are bounded, globally Lipschitz and 
\begin{equation}\label{eq2.1}
	\alpha_0 \leq \alpha(t) \leq \alpha_1 \quad \text { for some } \quad \alpha_1 \geq \alpha_0>0,
\end{equation}
\begin{equation}\label{eq2.2}
		\kappa_0 \leq \kappa(t) \leq \kappa_1 \quad \text { for some } \quad \kappa_1 \geq \kappa_0>0.
\end{equation}

We assume the following conditions for the growth and dissipation of the nonlinearity $f=\left(f_1, \ldots, f_N\right)$. We regard $f$ as a conservative vector field in $\mathbb{R}^N$, with each component function $f_i: \mathbb{R}^N \rightarrow \mathbb{R}$ being twice continuously differentiable for any $i \in{1, \ldots, N}$. Additionally, for some $\eta \in\left(0,\min\left\{\frac{\lambda_1(2\mu+\lambda)}{2},\lambda_1\right\} \right)$ and a positive constant $C_f$, where $\lambda_1>0$ is the first eigenvalue of the negative Laplacian operator with zero Dirichlet boundary conditions, we make the following assumptions
\begin{equation}\label{eq2.3}
	-C_f-\frac{\eta}{2}|u|^2 \leq \hat{f}(u):=\int_0^u f \mathrm{~d} \gamma \leq f(u) \cdot u+\frac{\eta}{2}|u|^2,
\end{equation}
with $\cdot$ denoting the standard dot product on $\mathbb{R}^N$, and $\int_0^u f \mathrm{~d} \gamma$ represents the line integral of $f$ along a piecewise smooth curve with initial point 0 and final point $u$, for any $u \in \mathbb{R}^N$, that is,
$$
\nabla \hat{f}(u)=f(u),
$$
where $\nabla \hat{f}$ stands for the gradient of $\hat{f}$ in the variables $u \in \mathbb{R}^N$. In addition, we shall assume throughout this paper that there exists a constant $C>0$ such that for every $i \in\{1, \ldots, n\}$ and $\xi=\left(\xi_1, \ldots, \xi_N\right) \in \mathbb{R}^N$,
\begin{equation}\label{eq2.4}
	\left|\nabla f_i(\xi)\right| \leqslant C\left(1+\sum_{i=1}^N\left|\xi_i\right|^{\rho-1}\right),
\end{equation}
and
\begin{equation}\label{eq2.5}
	\left|\partial_{x_i}^2 f_i(\xi)\right| \leqslant C,
\end{equation}
for some $1<\rho<\frac{N}{N-2}$, if $N \geq 3$.

The function $g(\cdot, s)$ is a given external force which is assumed translation bounded in $L_{l o c}^2\left(\mathbb{R} ; L^2(\Omega)\right)$, i.e., $g(\cdot) \in L_b^2\left(\mathbb{R} ; L^2(\Omega)\right)$,
\begin{equation}\label{eq2.6}
	\|g\|_{L_b^2\left(\mathbb{R} ; L^2(\Omega)\right)}^2=\sup _{t \in \mathbb{R}} \int_t^{t+1}\|g(s)\|_{L^2(\Omega)}^2 d s<\infty .
\end{equation} 

\textbf{\large {2.2 Well-posedness}}

In this subsection, based on \cite{FD1}, we show that the displacement can be decomposed into two parts: a curl free part and the divergence free part. To do this, we solve system over the space
$$
H_c=W_c \times E_c \times L^2(\Omega),
$$
and
$$
H_d=W_d \times E_d,
$$
where
$$
\begin{aligned}
	& W_c=\left\{u \in\left[H_0^1(\Omega)\right]^N ; \operatorname{curl} u=0\right\}, \\
	& E_c=\left\{u \in\left[L^2(\Omega)\right]^N ; \operatorname{curl} u=0\right\}, \\
	& W_d=\left\{u \in\left[H_0^1(\Omega)\right]^N ; \operatorname{div} u=0\right\},
\end{aligned}
$$
and
$$
E_d=\left\{u \in\left[L^2(\Omega)\right]^N ; \operatorname{div} u=0\right\} .
$$

If $u \in\left[H_0^1(\Omega)\right]^N$, then there exists functions $u^c, u^d \in\left[H_0^1(\Omega)\right]^N$, with curl $u^c=0$ and div $u^c=0$, such that
$$
u=u^c+u^d .
$$

So we have that
$$
\left[H_0^1(\Omega)\right]^N=W_c \oplus W_d, \quad W_c \cap W_d=\{0\}.
$$

Under the notations and conditions given above, we can then decompose problem (\ref{eq1.1})-(\ref{eq1.2}) into the following problems.

$\mathbf{Case\,\, curl \,\,u^c=0:}$
\begin{equation}\label{eq2.7}
	\left\{\begin{array}{l}
		\partial_t^2 u^c-\Delta_e u^c+\alpha(t) \nabla \theta+f(u)=\nabla p, \\
		\partial_t \theta-\kappa(t) \Delta \theta+\alpha(t) \operatorname{div} \partial_t u^c=g(x,t),
	\end{array} \quad \Omega \times[\tau,+\infty),\right.
\end{equation}
with initial conditions and boundary conditions
\begin{equation}\label{eq2.8}
	\left\{\begin{array}{l}
		u^c(x, \tau)=u_\tau^c(x),\, u_t^c(x, \tau)=v_\tau^c(x),\, \theta(x, \tau)=\theta_\tau(x),\quad x\in\Omega,\,\tau\in\mathbb{R}, \\
		u^c(x,t)=0,\,\theta(x, t)=0,\,x\in\partial\Omega,\,t\geq\tau,
	\end{array}\right.
\end{equation}
where $\Delta p=0$, $\nabla p$ and $u_t$ are orthogonal.

Denoting $U = (u^c,z^c,\theta)$ be the state vector with $v^c = u_t^c$, we can rewrite the system (\ref{eq2.7})-(\ref{eq2.8}) as an equivalent Cauchy problem in the product space $H_c$
\begin{equation}\label{eq2.9}
	\left\{\begin{array}{l}
		\frac{\mathrm{d} U}{\mathrm{~d} t}+A_c(t) U=F(U), t>\tau, \\
		U(\tau)=U_\tau,
	\end{array}\right.
\end{equation}
where $U=\left[\begin{array}{c}u^c \\ v^c \\ \theta\end{array}\right], U_\tau=\left[\begin{array}{c}u_\tau^c \\ v_\tau^c \\ \theta_\tau\end{array}\right]$, the linear operator $A_c(t): D\left(A_c\right) \subset H_c \rightarrow H_c$ is a unbounded linear operator defined by
$$
A_c(t)=\left[\begin{array}{ccc}
	0 & -I & 0 \\
	-\Delta_e & 0 & \alpha(t) \nabla \\
	0 & \alpha(t) \operatorname{div} & -\kappa(t) \Delta
\end{array}\right],
$$
with domain
$$
D\left(A_c(t)\right)=\left(W_c \cap\left[H^2(\Omega)\right]^N\right) \times W_c \times H_0^1(\Omega) \cap H^2(\Omega),
$$
and $F: H_c \rightarrow H_c$ is a nonlinear operator given by
$$
F(U)=\left[\begin{array}{c}
	0 \\
	-f^e\left(u^c\right) \\
	g(x,t)
\end{array}\right],
$$
where $f^e$ denotes the Nemyts\v{k}i operator associated with $f$; that is, for any $x \in \Omega$ and $t \geq \tau$
$$
f^e\left(u^c\right)(x, t)=f(u(x, t))=\left(f_1\left(u^c(x, t)\right), \ldots, f_N\left(u^c(x, t)\right)\right) .
$$

We can obtain these lemmas in \cite{FD2}.
\begin{Lemma}(\cite{FD2})\label{le2.1}
	 If the functions $f_i$ satisfy (\ref{eq2.4}), then there exists a constant $C>0$ such that for every $i=1, \ldots, n$, and $u=\left(u_1, \ldots, u_n\right), y=\left(y_1, \ldots, y_n\right) \in \mathbb{R}^N$, we have
	$$
	\left|f_i(u)-f_i(y)\right| \leqslant 2^{\rho-1} N|u-y|\left(1+\sum_{i=1}^N\left|u_i\right|^{\rho-1}+\sum_{i=1}^N\left|y_i\right|^{\rho-1}\right) .
	$$
	Consequently, there exists a constant $\bar{C}>0$ for any $U_1=\left(u_1, z_1, \theta_1\right), U_2=\left(u_2, z_2, \theta_2\right) \in H_c$ with $u_i=\left(u_{i 1}, \ldots, u_{i N}\right)$ and we deduce that
	$$
	\left\|F\left(U_1\right)-F\left(U_2\right)\right\|_{\left[H_0^1(\Omega)\right]^N \times\left[L^2(\Omega)\right]^N \times L^2(\Omega)} \leqslant \bar{C}\left\|U_1-U_2\right\|_{\left(H^1(\Omega)\right)^N}\left(1+\sum_{i=1}^2 \sum_{j=1}^N\left\|u_{i j}\right\|_{H^1(\Omega)}^{\rho-1}\right) .
	$$
\end{Lemma}
\begin{Lemma}(\cite{FD2})\label{le2.2}
	 If the functions $f_i$ satisfy (1.7), then the Nemyts\v{k}ii operators associated with $f_i, f_i^e: W_c \rightarrow$ $L^2(\Omega)$ are continuously differentiable and the derivative operators $D f_i^e:\left[H_0^1(\Omega)\right]^N \rightarrow \mathcal{L}\left(\left[H_0^1(\Omega)\right]^N, L^2(\Omega)\right)$ are Lipschitz continuous (in bounded subsets of $W_c$ ), for $N=3,4$. For $N>4$, there exists a constant $\eta \in(0,1)$ such that for every $u, v \in W_c$,
	$$
	\left\|D f_i^e(u)-D f_i^e(v)\right\|_{\mathcal{L}\left(\left[H_0^1(\Omega)\right]^N, L^2(\Omega)\right)} \leqslant c\|u-v\|_{W_c}^\eta .
	$$	
\end{Lemma}
 We define the energy functional $E(t)$ associated to system (\ref{eq1.1})-(\ref{eq1.2}) along a weak solution $U = (u,u_t,\theta)$ by
 \begin{equation}\label{eq2.10}
 	E(t)=\frac{1}{2}\|u_t\|^2 +\frac{1}{2}\mu\|\nabla u\|^2+\frac{1}{2}(\mu+\lambda)\|\operatorname{div} u\|^2+\frac{1}{2}\|\theta\|^2+ \int_{\Omega} \hat{f}(u) \mathrm{~d} x,
 \end{equation}
 and
  \begin{equation}\label{eq2.11}
 	E_c(t)=\frac{1}{2}\|u_t\|^2+\frac{1}{2}(2\mu+\lambda)\|\operatorname{div} u\|^2+\frac{1}{2}\|\theta\|^2+\int_{\Omega} \hat{f}(u) \mathrm{~d} x .
 \end{equation}

 If $U =(u,u_t,\theta)$ is a strong solution, then
 \begin{equation}\label{eq2.12}
 	\frac{d}{dt}E(t)=-\kappa(t)\|\nabla\theta\|^2+\left<g(t),\theta\right>.
 \end{equation}

To enhance the clarity of our results, we establish certain definitions. Let us define the notation $(\cdot, \cdot)_{H_0^1(\Omega)}$ to represent the inner product space in $H_0^1(\Omega)$ and denote by
$$
(\nabla u, \nabla v)=\sum_{i=1}^N\left(\nabla u_i, \nabla v_i\right)=\int_{\Omega} \nabla u \nabla v d x,
$$
the inner product in $\left(H_0^1(\Omega)\right)^N$. Indeed, to deal with problem (\ref{eq1.1})-(\ref{eq1.2}), we consider the Hilbert space $\left(H_0^1(\Omega)\right)^N$ equipped with the inner product
\begin{equation}\label{eq2.13}
	\left(u_1, u_2\right)_{\left(H_0^1(\Omega)\right)^N}=\int_{\Omega}\left[\mu \nabla u_1 \nabla u_2+(\lambda+\mu) \operatorname{div} u_1 \operatorname{div} u_2\right] d x .
\end{equation}
The norm induced by this inner product is equivalent to the usual one of $\left(H_0^1(\Omega)\right)^N$ and satisfies
\begin{equation}\label{eq2.14}
	\mu\|\nabla u\|_{\left(L^2(\Omega)\right)^N}^2 \leqslant\|u\|_{\left(H_0^1(\Omega)\right)^N}^2 \leqslant \tilde{c}\|\nabla u\|_{\left(L^2(\Omega)\right)^N}^2, \quad u \in\left(H_0^1(\Omega)\right)^N,
\end{equation}
where $\tilde{c}=\mu+N(\lambda+\mu)$.
\begin{Lemma}\label{le2.3}
	There exist constants $\beta_0$, $C_f>0$, such that
	\begin{equation}\label{eq2.15}
		E(t)\geq \beta_0\|(u,u_t,\theta\|^2_{H_c}-C_f|\Omega|.
	\end{equation}
\rm\textbf{Proof.} Let us consider the continuous functional $E_c: H_c \rightarrow \mathbb{R}$ defined using (\ref{eq2.10}) as $E_c(t) \approx E(t)$. Using (\ref{eq2.3}), we obtain
\begin{equation}
	\begin{aligned}
		E(t)& \geq \frac{1}{2}\left(\|u\|_{\left[H_0^1(\Omega)\right]^N}^2+\|z\|_{\left[L^2(\Omega)\right]^N}^2+\|\theta\|_{L^2(\Omega)}^2\right)-\frac{\eta}{2}\|u\|_{\left[L^2(\Omega)\right]^N}^2-C_\eta|\Omega|,\\
		&\geq \frac{1}{2}\left(\|u\|_{\left[H_0^1(\Omega)\right]^N}^2+\|z\|_{\left[L^2(\Omega)\right]^N}^2+\|\theta\|_{L^2(\Omega)}^2\right)-\frac{\eta}{2\lambda_1}\|u\|_{\left[H_0^1(\Omega)\right]^N}^2-C_\eta|\Omega|.
	\end{aligned}
\end{equation}
Taking $\beta_0=\frac{1}{2}(1-\frac{\eta}{\lambda_1})$, we can easily get (\ref{eq2.15}).$\hfill\qedsymbol$
\end{Lemma}
From Lemma \ref{le2.1}, Lemma \ref{le2.2} and \cite{FD1}, we know the problem (\ref{eq2.7})-(\ref{eq2.8}) has a unique local solution $U(t)$ in $H_c$ satisfying
the initial condition $U(\tau)=U_{\tau} \in H_c$. Now, suppose $t_{\max }<+\infty$, by (\ref{eq2.12}) and Young's inequality\cite{YQ25,YQ26}, we get
$$
E(t) \leq E(\tau)+\frac{1}{2 k_0} \int_\tau^{t_{\max }}\|g(s)\|^2 d s, \quad \forall t \in\left[\tau, t_{\max }\right) .
$$
By (\ref{eq2.15}) in Lemma \ref{le2.3}, we have

$$
\left\|\left(u(t), u_t(t),\theta\right)\right\|_{{H_c}}^2 \leq \frac{1}{\beta_0}\left(E(\tau)+\int_\tau^{t_{\max }}\|g(s)\|^2 d s+C_f|\Omega|\right)=C_\tau,
$$
for all $t \in\left[\tau, t_{\max }\right)$. Therefore, $t_{\max }=+\infty$ and $U$ is a global solution.

$\mathbf{Case\,\,div\,\,u^d=0:}$
\begin{equation}
	 u^d_{tt}-\mu \Delta u^d=\nabla r, \quad \Omega \times[\tau,+\infty),
\end{equation}
with initial conditions and boundary conditions
\begin{equation}
	\left\{\begin{array}{l}
		 u^d(x, \tau)=u_\tau^d(x), u^d_{t}(x, \tau)=v_\tau^d(x),\\
		u^d(x, t)=0, x \in \partial \Omega, t \geq \tau,
	\end{array}\right.
\end{equation}
where $\Delta r=0$.

It should be noted that this case is identical to that in \cite{FD1}, hence we omit it. 

Combining the above analysis, the following theorem may be established by a straightforward extension of \cite{FD1}, so we omit its proof.
\begin{Theorem}(\cite{FD1})\label{th2.4}
	Under assumption (\ref{eq2.1}-(\ref{eq2.6}), the following statements hold.
	
	If $\left(u_\tau, v_\tau, \theta_\tau\right) \in H_c$, then problem (\ref{eq1.1})-(\ref{eq1.2}) has a unique mild solution $\left(u, \partial_t u, \theta\right) \in C\left([\tau, \infty), H_c\right)$;
	
	 Moreover, if $\left(u_\tau, v_\tau, \theta_\tau\right) \in D\left(A_c(t)\right)$, then the above mild solution $U$ is a classical solution and satisfies	
	$$
	\left(u, \partial_t u, v\right) \in C\left([\tau, \infty),\left(W_c \cap\left[H^2(\Omega)\right]^N\right) \times W_c \times H_0^1(\Omega) \cap H^2(\Omega)\right) \cap C^1\left([\tau, \infty), H_c\right) .
	$$	
\end{Theorem}

\section{ {Uniformly (w.r.t. $G \in \Sigma$ ) absorbing set in $H_c$}}

 \quad In this section, we focus on demonstrating the exisence of a uniform absorbing set within the evolution process initiated by the system of equations (\ref{eq2.7})-(\ref{eq2.8}) in the space $H_c$. To set the stage, we begin by revisiting some basic concepts about non-autonomous systems, with \cite{VV3} serving as a comprehensive reference for further exploration.
 
\textbf{\large {3.1 Abstract results}}
\begin{Definition}(\cite{VV3})\label{de3.1}
	Let $X$ be a Banach space, and $\Sigma$ a parameter set. The operators $\left\{U_\sigma(t, \tau)\right\}, \sigma \in \Sigma$, are said to be a family of processes in $X$ with symbol space $\Sigma$ if for any $\sigma \in \Sigma$
	$$
	\begin{aligned}
		& U_\sigma(t, s) \circ U_\sigma(s, \tau)=U_\sigma(t, \tau), \quad \forall t \geqslant s \geqslant \tau, \tau \in \mathbb{R}, \\
		& U_\sigma(\tau, \tau)=\operatorname{Id} \text { (identity), } \quad \forall \tau \in \mathbb{R} .
	\end{aligned}
	$$
	Let $\{T(s)\}_s \geqslant 0$ be the translation semigroup on $\Sigma$. We say that a family of processes $\left\{U_\sigma(t, \tau)\right\}, \sigma \in \Sigma$, satisfies the translation identity if
	$$
	\begin{aligned}
		& U_\sigma(t+s, \tau+s)=U_{T(s) \sigma}(t, \tau), \quad \forall \sigma \in \Sigma, t \geqslant \tau, \tau \in \mathbb{R}, s \geqslant 0, \\
		& T(s) \Sigma=\Sigma, \quad \forall s \geqslant 0.
	\end{aligned}
	$$
\end{Definition}
\begin{Definition}(\cite{VV3})\label{de3.2}
	A bounded set $B_0 \in \mathcal{B}(X)$ is said to be a bounded uniformly (w.r.t. $\sigma \in \Sigma$ ) absorbing set for $\left\{U_\sigma(t, \tau)\right\}, \sigma \in \Sigma$, if for any $\tau \in \mathbb{R}$ and $B \in \mathcal{B}(X)$ there exists $T_0=T_0(B, \tau)$ such that $\bigcup\limits_{\sigma \in \Sigma} U_\sigma(t ; \tau) B \subset B_0$ for all $t \geqslant T_0$.
\end{Definition}
\begin{Definition}(\cite{VV3})\label{de3.3}
	A set $\mathcal{A} \subset X$ is said to be uniformly (w.r.t $\sigma \in \Sigma$ ) attracting for the family of processes $\left\{U_\sigma(t, \tau)\right\}, \sigma \in \Sigma$, if for any fixed $\tau \in \mathbb{R}$ and any $B \in \mathcal{B}(X)$
	$$
	\lim _{t \rightarrow+\infty}\left(\sup _{\sigma \in \Sigma} \operatorname{dist}\left(U_\sigma(t ; \tau) B ; \mathcal{A}\right)\right)=0,
	$$
	where $\operatorname{dist}(\cdot, \cdot)$ is the usual Hausdorff semidistance in $X$ between two sets.
	In particular, a closed uniformly attracting set $\mathcal{A}_{\Sigma}$ is said to be the uniform (w.r.t. $\sigma \in \Sigma$ ) attractor of the family of processes $\left\{U_\sigma(t, \tau)\right\}, \sigma \in \Sigma$, if it is contained in any closed uniformly attracting set (minimality property).
\end{Definition}
 \begin{Theorem}(\cite{IC4})\label{th3.4}
 	 Let $\Sigma$ be defined as before and $G_0 \in \mathbb{E}$, then
 	
 	(1) $G_0$ is a translation compact in $\hat{\mathbb{E}}$ and for any $G \in \mathbb{E}=H\left(G_0\right)$ is also a translation compact in $\hat{\mathbb{E}}$, moreover, $H(G) \subseteq$ $H\left(G_0\right)$;
 	
 	(2) The set $H\left(G_0\right)$ is bounded in $L_{\text {loc }}^2\left(R^{+},H_c\right)$, such that	
 	$$
 	\eta_G(h) \leq \eta_{G_0}(h)<+\infty, \text { for any } G \in \Sigma .
 	$$
 \end{Theorem}
 
\textbf{\large {3.2 Existence of uniformly (w.r.t. $G \in \Sigma$ ) absorbing set in $H_c$}} 

We will now present a theorem regarding the existence of uniformly absorbing sets. Before stating the main result, we prove several important lemmas as follows.
\begin{Lemma}\label{le3.5}
	Let $(u,\theta)$ be the solution to (\ref{eq2.7})-(\ref{eq2.8}), the energy function defined by (\ref{eq2.11}), then 
	\begin{equation}\label{eq3.1}
		\frac{d}{dt}E(t)\leq -\frac{k_0}{2}\|\nabla\theta\|^2+\frac{1}{2k_0\lambda_1}\|g\|^2.	 
	\end{equation}
    \rm\textbf{Proof.} From (\ref{eq2.12}), we know
    $$
    	\frac{d}{dt}E(t)=-\kappa(t)\|\nabla\theta\|^2+\left<g(t),\theta\right>.
   $$
   
   By applying Young's inequality, we can identify a positive constant $\epsilon_1$ such that the following inequality holds,
   $$
        \int_{\Omega}g\cdot \theta dx\leq \frac{\epsilon_1}{2}\|\theta\|^2+\frac{1}{2\epsilon_1}\|g\|^2.
   $$
   
   Utilizing the Poinc\'{a}re inequality, assumption (\ref{eq2.2}) and setting $\epsilon_1 = k_0 \lambda_1$, we can readily deduce inequality (\ref{eq3.1}).$\hfill\qedsymbol$
\end{Lemma}
\begin{Lemma}\label{le3.6}
	 	Let $(u,\theta)$ be the solution to (\ref{eq2.7})-(\ref{eq2.8}), the function $F_1(t)$ defined by
	 	\begin{equation}\label{eq3.2}
	 		F_1(t)=\int_{\Omega}u\cdot u_tdx,
	 	\end{equation}
 	 satisfies
 	 \begin{equation}\label{eq3.3}
 	 	\frac{d}{dt}F_1(t)\leq \|u_t\|^2-\frac{2\mu+\lambda}{2}\|\operatorname{div}u\|^2+\frac{2\lambda_1\alpha_0}{4\eta}\|\nabla\theta\|^2-\int_{\Omega} \hat{f}(u) \mathrm{~d} x.
 	 \end{equation}
  \rm\textbf{Proof.} Multiplying equation $(\ref{eq2.7})_1$ by $u$ and applying theorthogonality condition, we get
  \begin{equation}\label{eq3.4}
  	\frac{d}{dt}=\|u_t\|^2-\frac{2\mu+\lambda}{2}\|\operatorname{div}u\|^2-\alpha(t)\int_{\Omega}\nabla\cdot \theta \mathrm{~d}x-\int_{\Omega} f(u)\cdot u\mathrm{~d}x.
  \end{equation}

 Using assumption (\ref{eq2.1}), immersion  $H_0^1(\Omega) \hookrightarrow L^2(\Omega)$ and Young's inequality
 \begin{equation}\label{eq3.5}
 	\begin{aligned}
 		-\alpha(t)\int_{\Omega}\nabla\cdot \theta \mathrm{~d}x&\leq -\alpha_0\int_{\Omega}\nabla\theta\cdot u \mathrm{~d}x\leq \alpha_0\int_{\Omega}|\nabla\theta|\cdot |u| \mathrm{~d}x\\
 		&\leq \frac{2\lambda_1\alpha_0^2}{4\eta}\|\nabla\theta\|^2+\frac{\eta}{2\lambda_1}\|\operatorname{div}u\|^2.
 	\end{aligned}
 \end{equation} 

Using assumption (\ref{eq2.3}), we have
\begin{equation}\label{eq3.6}
	-\int_{\Omega} f(u)\cdot u\mathrm{~d}x\leq -\int_{\Omega}\hat{f}(u)\mathrm{~d}x+\frac{\eta}{2\lambda_1}\|\operatorname{div}u\|^2.
\end{equation}

 Inserting the inequality (\ref{eq3.5}) and (\ref{eq3.6}) in (\ref{eq3.4}), we readily obtain equation (\ref{eq3.3}). $\hfill\qedsymbol$
\end{Lemma}

Now, we introduce the multiplicators $\phi$ and $w$ given by
$$
-\Delta \phi=\operatorname{div} u, \quad-\Delta w=\theta, \quad \text { in } \quad \Omega,
$$
with
$$
\phi(x, t)=0, \quad w(x, t)=0, \quad \text { on } \quad \Gamma=\partial \Omega .
$$
\begin{Lemma}\label{le3.7}
	Under above notations, let $(u,\theta)$ be the solution to (\ref{eq2.7})-(\ref{eq2.8}), then 
	\begin{equation}\label{eq3.7}
		F_2(t)=\int_{\Omega} \theta \partial_t \phi \mathrm{~d} x,
	\end{equation}
satisfies
\begin{equation}\label{eq3.8}
	\frac{d}{dt}F_2(t)\leq -\frac{\alpha_0}{2k}\|u_t\|^2+C_{\epsilon,\delta}\|\nabla\theta\|^2+\epsilon\|\operatorname{div}\|^2+\delta\int_{\Gamma}|\operatorname{div} u|^2\mathrm{~d} \Gamma+\frac{1}{\alpha_0\lambda_1}\|g(t)\|^2.
\end{equation}
\rm\textbf{Proof.}
By applying integration by parts to equation $(\ref{eq2.7})_2$ and performing direct calculations, we can obtain the following result
\begin{equation}\label{eq3.9}
	\begin{aligned}
	\frac{\mathrm{d}}{\mathrm{~d} t}F_2(t) & =-\alpha(t) \int_{\Omega}\left|\nabla \partial_t \phi\right|^2 \mathrm{d} x-\kappa(t) \int_{\Omega} \nabla \theta \cdot \nabla \partial_t \phi \mathrm{~d} x-(2 \mu+\lambda) \int_{\Gamma} \frac{\partial w}{\partial \nu} \operatorname{div} u \mathrm{~d} \Gamma \\
		&-(2 \mu+\lambda) \int_{\Omega} \theta \operatorname{div} u \mathrm{~d} x+\alpha(t) \int_{\Omega} \nabla w \cdot \nabla \theta \mathrm{~d} x+\int_{\Omega} \nabla w \cdot f(u) \mathrm{d} x+\int_{\Omega}g(t)\partial_t\phi\mathrm{~d} x .
	\end{aligned}
\end{equation}

Now we estimate the terms on the right-hand side of equation (\ref{eq3.9}).

Applying assumption (\ref{eq2.1}), there exists a constant $k>0$, such that
\begin{equation}\label{eq3.10}
		-\alpha(t) \int_{\Omega}\left|\nabla \partial_t \phi\right| \mathrm{d} x\leq -\alpha_0\|\nabla\partial_t\phi\|^2\leq-\frac{\alpha_0}{k}\|u_t\|^2.
\end{equation}

Using assumption (\ref{eq2.2}) and Young's inequality, we have
\begin{equation}\label{eq3.11}
	\begin{aligned}
		-\kappa(t) \int_{\Omega} \nabla \theta \cdot \nabla \partial_t \phi \mathrm{~d} x\leq k_1\|\nabla \theta\|\| \nabla\partial_t\phi\|&\leq \frac{k_1^2}{\alpha_0}\|\nabla\theta\|^2+\frac{\alpha_0}{4}\|\nabla\partial_t\phi\|^2\\
		& \leq \frac{k_1^2}{\alpha_0}\|\nabla\theta\|^2+\frac{\alpha_0}{4k}\|u_t\|^2.
	\end{aligned}
\end{equation}

The following inequalities can be directly obtained from \cite{FD1}, and their proofs are omitted here
\begin{equation}\label{eq3.12}
	-(2 \mu+\lambda) \int_{\Gamma} \frac{\partial w}{\partial \nu} \operatorname{div} u \mathrm{~d} \Gamma \leq C_\delta\|\nabla \theta\|^2+\delta \int_{\Gamma}|\operatorname{div} u|^2 \mathrm{~d} \Gamma ,
\end{equation}
\begin{equation}\label{eq3.13}
	-(2 \mu+\lambda) \int_{\Omega} \theta \operatorname{div} u \mathrm{~d} x \leq C_\epsilon^{\prime}\|\nabla \theta\|^2+\frac{\epsilon}{2} \int_{\Omega}|\operatorname{div} u|^2 \mathrm{~d} x ,
\end{equation}
\begin{equation}\label{eq3.14}
	\alpha(t) \int_{\Omega} \nabla w \cdot \nabla \theta \mathrm{~d} x \leq \frac{\alpha_1}{\lambda_1}\|\nabla \theta\|^2 ,
\end{equation}
\begin{equation}\label{eq3.15}
	\int_{\Omega} \nabla w \cdot f(u) d x \leq C_\epsilon^{\prime \prime}\|\nabla \theta\|^2+\frac{\epsilon}{2} \int_{\Omega}|\operatorname{div} u|^2 \mathrm{~d} x,
\end{equation}
where $\epsilon$ and $\delta$ are positive constants that will be chosen later.

Utilizing the Young's inequality and Poinc\'{a}re inequality, we have
\begin{equation}\label{eq3.16}
	\begin{aligned}
		\int_{\Omega}g(t)\partial_t\phi\mathrm{~d} x&\leq \|g(t)\|\|\partial_t\phi\|
		\leq\lambda_1^\frac{-1}{2}\|g(t)\|\nabla\|\partial_t\phi\| \\&\leq\frac{1}{\alpha_0\lambda_1}\|g(t)\|^2+\frac{\alpha_0}{4}\|\partial_t\phi\|^2\leq\frac{1}{\alpha_0\lambda_1}\|g(t)\|^2+\frac{\alpha_0}{4k}\|u_t\|^2.
	\end{aligned}	
\end{equation}

Incorporating inequalities (\ref{eq3.10})- (\ref{eq3.16}) into equation (\ref{eq3.9}), we straightforwardly derive equation (\ref{eq3.8}).
$\hfill\qedsymbol$
\end{Lemma} 
\begin{Lemma}\label{le3.8}
	Let us denote by $q$ a $C^2$ function such that $q=\nu$ over $\Gamma$ and $(u,\theta)$ be the solution to (\ref{eq2.7})-(\ref{eq2.8}), then 
	\begin{equation}\label{eq3.17}
		F_3(t)==-\sum_{j=1}^N \sum_{i=1}^N\left\{\int_{\Omega} \partial_t u_j q_i \frac{\partial u_i}{\partial x_j} \mathrm{~d} x+\int_{\Omega} \partial_t u_j \frac{\partial q_i}{\partial x_j} u_i \mathrm{~d} x\right\} ,
	\end{equation}
satisfies
\begin{equation}\label{eq3.18}
	\frac{d}{dt}F_3(t)\leq \frac{3M_q}{2}\|u_t\|^2+P\|\operatorname{div} u\|^2+C_{\epsilon}^{'''}\|\nabla\theta\|^2-\frac{2\mu+\lambda}{2}\int_{\Gamma}|\operatorname{div} u|^2\mathrm{~d} \Gamma+M_2,
\end{equation}
where the constants $M_q$, $P$, $M_2>0$, will be defined  in the subsequent sections.

\rm\textbf{Proof.}
Multiplying equation $(\ref{eq2.7})_1$ by $\nabla(q\cdot u)$ and utilizing the orthogonality condition, we can arrive at the following equality through straightforward calculations
\begin{equation}\label{eq3.19}
	\begin{aligned}
		\frac{d}{d} F_3(t)= & -\sum_{j=1}^N \sum_{i=1}^N \int_{\Omega} \partial_t u_j \frac{\partial q_i}{\partial x_j}\left(\partial_t u_i\right) \mathrm{d} x+\frac{1}{2} \sum_{i=1}^N \int_{\Omega} \operatorname{div} q\left[\left|\partial_t u_i\right|^2-(2 \mu+\lambda)|\operatorname{div} u|^2\right] \mathrm{d} x \\
		& -\frac{(2 \mu+\lambda)}{2} \int_{\Gamma}|\operatorname{div} u|^2 \mathrm{~d} \Gamma+(2 \mu+\lambda) \sum_{j=1}^N \int_{\Omega} \operatorname{div} u \sum_{i=1}^N\left[\frac{\partial^2 q_i}{\partial x_j^2} u_i+2 \frac{\partial q_i}{\partial x_j} \frac{\partial u_i}{\partial x_j}\right] \mathrm{d} x \\
		& +\int_{\Omega} \alpha(t) \nabla \theta \cdot \nabla(q \cdot u) \mathrm{d} x-\int_{\Omega} f(u) \cdot \nabla(q \cdot u) \mathrm{d} x .
	\end{aligned}
\end{equation}

Next, we will estimate the terms on the right-hand side of equation (\ref{eq3.19}). We first denote 
$$M_q\geq \max\limits_{x\in\bar{\Omega}}\{\|Dq\|,\|\operatorname{div} q\|,\|\nabla q\|\}.
$$

Applying Young's inequality, we get
\begin{equation}\label{eq3.20}
	-\sum_{j=1}^N \sum_{i=1}^N \int_{\Omega} \partial_t u_j \frac{\partial q_i}{\partial x_j}\left(\partial_t u_i\right) d x=-\int_{\Omega} \partial_t u \cdot(D q)^{\top} \partial_t u d x \leq M_q\left\|u_t\right\|^2.
\end{equation}

Based on the definition of $M_q$, we derive
\begin{equation}\label{eq3.21}
	\begin{aligned}
		\frac{1}{2} \sum_{i=1}^N \int_{\Omega} \operatorname{div} q\left[\left|\partial_t u_i\right|^2-(2 \mu+\lambda)|\operatorname{div} u|^2\right] \mathrm{d} x
		&=\frac{1}{2}\int_{\Omega}\operatorname{div}q\left[|\partial_t u|^2-(2 \mu+\lambda)|\operatorname{div} u|^2 \right]\mathrm{d} x\\
		&\leq \frac{M_q}{2}\|u_t\|^2-\frac{M_q(2 \mu+\lambda)}{2}\|\operatorname{div} u\|^2,
	\end{aligned}
\end{equation}
and
\begin{equation}\label{eq3.22}
	\begin{aligned}
		&(2 \mu+\lambda) \sum_{j=1}^N \int_{\Omega} \operatorname{div} u \sum_{i=1}^N\left[\frac{\partial^2 q_i}{\partial x_j^2} u_i+2 \frac{\partial q_i}{\partial x_j} \frac{\partial u_i}{\partial x_j}\right] \mathrm{d} x\\
		&\quad\quad=(2 u+\lambda) \int_{\Omega} \operatorname{div} u \sum_{i=1}^N(\Delta q_i u_i+2 \nabla q_i \cdot \nabla u_i) \mathrm{d}x\\
		&\quad\quad\leq 3(2 u+\lambda) M_q\|\operatorname{div} u \|^2.
	\end{aligned}
\end{equation}

Applying Young's inequality again, we obtain
\begin{equation}\label{eq3.23}
	\begin{aligned}
		\int_{\Omega} \alpha(t) \nabla \theta \cdot \nabla(q \cdot u) \mathrm{d} x&=\int_\Omega \alpha(t) \cdot \nabla \theta \cdot[(\nabla q) \cdot u+q \cdot \nabla u] d x \\
		&\leq C_\epsilon^{\prime \prime \prime}\|\nabla \theta\|^2+\epsilon\|\operatorname{div} u\|^2,
	\end{aligned}
\end{equation}
and
\begin{equation}\label{eq3.24}
	\begin{aligned}
		-\int_{\Omega} f(u) \cdot \nabla(q \cdot u) \mathrm{d} x&=-\int_{\Omega} f(u) \cdot(\nabla q \cdot u) \mathrm{d} x-\int_{\Omega} f(u) \cdot(q \cdot \nabla u)  \mathrm{d}x \\
		&\leq \left|-\int_{\Omega} f(u) \cdot(\nabla q \cdot u) \mathrm{d} x\right|+\left|-\int_{\Omega} f(u) \cdot(q \cdot \nabla u)  \mathrm{d}x\right|\\
		&\leq M_q \int_{\Omega}|f(u)||u|\mathrm{d} x+M_q \int_{\Omega}|f(u)||\nabla u|\mathrm{d} x\\
		&\leq \frac{1}{2}\int_{\Omega}|f(u)|^2\mathrm{d} x+\frac{M_q^2}{2\lambda_1}\|\operatorname{div} u\|^2+\frac{1}{2}\int_{\Omega}|f(u)|^2\mathrm{d} x+\frac{M_q^2}{2}\|\operatorname{div} u\|^2\\
		&=\int_{\Omega}|f(u)|^2\mathrm{d} x+\frac{M_q^2
		+M_q^2\lambda_1}{2\lambda_1}\|\operatorname{div} u\|^2.
	\end{aligned}
\end{equation}

Thanks to auusmption (\ref{eq2.4}), there exists a constant $M_1>0$ such that
$$
\int_{\Omega}|f(u)|^2 d x \leqslant M_1 \int_{\Omega}|u|^2 d x+M_1 \sum_{i=1}^N \int_{\Omega}\left|u_i\right|^{2 \rho} d x .
$$

Since $1<\rho<\frac{N}{N-2}$, it follows that  $H^1(\Omega) \hookrightarrow L^{2 \rho}(\Omega)$. Consequently, we obtain the following inequality
\begin{equation}\label{eq3.25}
	\int_{\Omega}|f(u)|^2 d x \leqslant M_1 \int_{\Omega}|u|^2 d x+\bar{M_1} \leqslant \frac{\bar{M_1}}{\lambda_1} \int_{\Omega}|\nabla u|^2 d x+M_2,
\end{equation}
whenever $\|u\|_{\left(H^1(\Omega)\right)^N} \leqslant r$ (as in \cite{AN5}).

Substituting equation (\ref{eq3.25}) into equation (\ref{eq3.24}), we can derive
\begin{equation}\label{eq3.26}
	\left|-\int_{\Omega} f(u) \cdot \nabla(q \cdot u) \mathrm{d} x\right|\leq \frac{2\bar{M_1}+M_q^2+\lambda_1M_q^2}{2\lambda_1}\|\operatorname{div} u\|^2+M_2.
\end{equation}

Finally, inserting equations (\ref{eq3.20})-(\ref{eq3.23}) and (\ref{eq3.26}) into equation (\ref{eq3.19}), and taking $P=\frac{5(2\mu+\lambda)M_q}{2}+\epsilon+\frac{2\bar{M_1}+M_q^2+\lambda_1M_q^2}{2\lambda_1}$, we obtain equation (\ref{eq3.18}).
$\hfill\qedsymbol$
\end{Lemma} 
 
 Now, we aim to construct a Lyapunov functional, which will be instrumental in achieving the desired result.
 \begin{Lemma}\label{le3.9}
 	For some constants $N_0,N_1,N_2,N_3 > 0$, we define a Lyapunov function $\mathcal{L}(t)$ as
 	\begin{equation}\label{eq3.27}
 		\mathcal{L}(t)=N_0E(t)+N_1F_1(t)+N_2F_2(t)+N_3F_3(t),
 	\end{equation}
satisfying
\begin{equation}\label{eq3.28}
	\frac{d}{dt}\mathcal{L}(t)\leq -\xi E(t)+\tilde{C}\|g\|^2+\tilde{M},
\end{equation} 
where $\xi$, $\tilde{M}$ and $\tilde{C}$ are positive constants that will be determined at a later stage. 

\rm\textbf{Proof.} 
By (\ref{eq3.1}), (\ref{eq3.2}), (\ref{eq3.7}) and  (\ref{eq3.18}), we can get
\begin{equation}\label{eq3.29}
	\begin{aligned}
		\frac{d}{dt}\mathcal{L}(t)&\leq N_0\left(-\frac{k_0}{2}\|\nabla\theta\|^2+\frac{1}{2k_0\lambda_1}\|g\|^2\right)\\
		&\quad +N_1\left(\|u_t\|^2-\frac{2\mu+\lambda}{2}\|\operatorname{div}u\|^2+\frac{2\lambda_1\alpha_0}{4\eta}\|\nabla\theta\|^2-\int_{\Omega} \hat{f}(u) \mathrm{~d} x\right)\\
		&\quad +N_2\left(-\frac{\alpha_0}{2k}\|u_t\|^2+C_{\epsilon,\delta}\|\nabla\theta\|^2+\epsilon\|\operatorname{div}\|^2+\delta\int_{\Gamma}|\operatorname{div} u|^2\mathrm{~d} \Gamma+\frac{1}{\alpha_0\lambda_1}\|g(t)\|^2\right)\\
		&\quad+N_3\left(\frac{3M_q}{2}\|u_t\|^2+P\|\operatorname{div}u\|^2+C_{\epsilon}^{'''}\|\nabla\theta\|^2-\frac{2\mu+\lambda}{2}\int_{\Gamma}|\operatorname{div} u|^2\mathrm{~d} \Gamma+M_2\right)\\
		&=-\left(\frac{k_0N_0}{2}-\frac{2\lambda_1\alpha_0N_1}{4\eta}-N_2C_{\epsilon,\delta}-N_3C_{\epsilon}^{'''}\right)\|\nabla\theta\|^2 -\left(\frac{\alpha_0N_2}{2k}-N_1-\frac{3M_qN_3}{2}\right)\|u_t\|^2\\
		&\quad -\left(\frac{N_1(2\mu+\lambda)}{2}-N_2\epsilon-PN_3\right)\|\operatorname{div}u\|^2-\left(\frac{N_3(2\mu+\lambda)}{2}-N_2\delta\right)\int_{\Gamma}|\operatorname{div} u|^2\mathrm{~d} \Gamma\\
		&\quad+\left(\frac{N_0}{2k_0\lambda_1}+\frac{N_2}{\alpha_0\lambda_1}\right)\|g(t)\|^2+N_3M_2.  
	\end{aligned}
\end{equation}

We deduce that if the following inequalities hold{\tiny }
\begin{equation}\label{eq3.30}
	\left\{\begin{array}{l}
		\frac{k_0N_0}{2}-\frac{2\lambda_1\alpha_0N_1}{4\eta}-N_2C_{\epsilon,\delta}-N_3C_{\epsilon}^{'''} \geq 0, \\
		\frac{\alpha_0N_2}{2k}-N_1-\frac{3M_qN_3}{2}\geq 0, \\
		\frac{N_1(2\mu+\lambda)}{2}-N_2\epsilon-PN_3\geq 0, \\
	    \frac{N_3(2\mu+\lambda)}{2}-N_2\delta\geq 0.
	\end{array}\right.
\end{equation}

We first denote $N_3=1$, $N_2=\frac{2\mu+\lambda}{3\delta}$, $N_1=\frac{
(2\mu+\lambda)\alpha_0}{6k\delta}-2M_q$,
then
$$
	\frac{N_1(2\mu+\lambda)}{2}-N_2\epsilon-PN_3=-\frac{(2\mu+\lambda)^2\alpha_0}{12k\delta}+(2\mu+\lambda)M_q+\frac{(2\mu+\lambda)\epsilon}{3\delta}+P.
$$

Let us denote $\delta=\frac{(2\mu+\lambda)^2\alpha_0}{12k[(2\mu+\lambda)M_q+2+P]}$, $\epsilon=\frac{3\delta}{2\mu+\lambda}$, which implies
$$
	N_1=\frac{2[(2\mu+\lambda)M_q+2+P]}{(2\mu+\lambda)}-2M_q=\frac{4+2P}{2\mu+\lambda}>0,
$$
$$
	N_2=\frac{4k[(2\mu+\lambda)M_q+2+P]}{(2\mu+\lambda)\alpha_0}=\frac{2\mu+\lambda}{3\delta}>0.
$$

Inserting $N_1$, $N_2$ in (\ref{eq3.30}), we can obtain
\begin{equation}\label{eq3.31}
	\left\{\begin{array}{l}
		\frac{\alpha_0N_2}{2k}-N_1-\frac{3M_qN_3}{2}=\frac{M_q}{2}\geq 0, \\
		\frac{N_1(2\mu+\lambda)}{2}-N_2\epsilon-PN_3=1, \\
		\frac{N_3(2\mu+\lambda)}{2}-N_2\delta=\frac{2\mu+\lambda}{6}\geq 0.
	\end{array}\right.
\end{equation}

Due to the embedding $H_0^1(\Omega) \hookrightarrow L^2(\Omega)$ with $\|\theta\|^2 \leq \frac{1}{\lambda_1}\|\nabla \theta\|^2$, we can choose $N_0$ sufficiently large such that
\begin{equation}\label{eq3.32}
	\lambda_1\left[\frac{k_0N_0}{2}-\frac{2\lambda_1\alpha_0N_1}{4\eta}-N_2C_{\epsilon,\delta}-N_3C_{\epsilon}^{'''}\right]\geq \max\left\{\frac{M_q}{2},1,\frac{2\mu+\lambda}{6}\right\}.
\end{equation}

From equations (\ref{eq3.31}) and (\ref{eq3.32}), we have verified that equation (\ref{eq3.30}) indeed holds true. 

In conclusion, we define  $\xi=\min\left\{\frac{M_q}{2},1,\frac{2\mu+\lambda}{6}\right\}$, $\tilde{C}=\frac{N_0}{2k_0\lambda_1}+\frac{N_2}{\alpha_0\lambda_1}$, $\tilde{M}=N_3M_2$. With these definitions, we can immediately derive equation (\ref{eq3.28}).
$\hfill\qedsymbol$
 \end{Lemma}
 
 In the sequel, we will establish the existence of the uniform (w.r.t. $G \in \Sigma$ ) absorbing set in $H_c$. Let $R_\tau=[\tau,+\infty), \tau \geq 0$, we will consider the following system
 
 \begin{equation}\label{eq3.33}
 	\left\{\begin{array}{l}
 		u_t-v=0,\\
 		\partial_t^2 u^c-\Delta_e u^c+\alpha(t) \nabla \theta+f(u)-\nabla p=0, \\
 		\partial_t \theta-\kappa(t) \Delta \theta+\alpha(t) \operatorname{div} \partial_t u^c=g(x,t),
 	\end{array} \quad \Omega \times[\tau,+\infty),\right.
 \end{equation}
 with initial and boundary conditions
 \begin{equation}\label{eq3.34}
 	\left\{\begin{array}{l}
 		u^c(x, \tau)=u_\tau^c(x),\, u_t^c(x, \tau)=v_\tau^c(x),\, \theta(x, \tau)=\theta_\tau(x),\quad x\in\Omega,\,\tau\in\mathbb{R}, \\
 		u^c(x,t)=0,\,\theta(x, t)=0,\,x\in\partial\Omega,\,t\geq\tau.
 	\end{array}\right.
 \end{equation}
 
 Let
 \begin{equation}\label{eq3.35}
 	 G=(0,0,g)^T \in \mathbb{E}=L^2\left(R_\tau,H_c\right) .
 \end{equation}
 
 For any $\left(u_\tau, v_{ \tau}, \theta_\tau\right) \in H_c$ and any $G \in \mathbb{E}, t \geq \tau, \tau \geq 0$, we define
 $$
 U_G(t, \tau):\left(u_\tau, v_{ \tau},  \theta_\tau\right) \in H_c \rightarrow\left(u(t), u_t(t), \theta(t)\right)=U_G(t, \tau)\left(u_\tau, v_{ \tau},  \theta_\tau\right),
 $$
 where $\left(u(t), u_t(t), \theta(t)\right)$ solves (\ref{eq3.33})-(\ref{eq3.34}).
 Define the hull of $F_0 \in \mathbb{E}$ as
 $$
 \Sigma=H\left(G_0\right)=\left[G_0(t+h) \mid h \in R^{+}\right]_{\mathbb{E}},
 $$
 where $[\cdot]_{\mathbb{E}}$ denotes the closure in Banach space $\mathbb{E}$.
 Note that
 $$
 G_0=(0,0,g_0) \in \mathbb{E} \subseteq \hat{\mathbb{E}},
 $$
 from assumption (\ref{eq2.6}), we deduce that $G_0$ is a translation compact function in $\hat{\mathbb{E}}$ in the weak topology, implying that $H\left(G_0\right)$ is compact in $\hat{\mathbb{E}}$. Subsequently, we consider the Banach space $L_{l o c}^2\left(R^{+}, \mathbb{E}_1\right)$, which consists of functions $\sigma(s), s \in R^{+}$, taking values in the Banach space $\mathbb{E}_1$ that are locally $p$-power integrable in the Bochner sense. In particular, for any $\left[t_1, t_2\right] \subseteq R^{+}$,
 $$
 \int_{t_1}^{t_2}\|\sigma(s)\|_{\mathbb{E}_1}^p d s<+\infty.
 $$
 
 Let $\sigma(s) \in L_{l o c}^2\left(R^{+}, \mathbb{E}_1\right)$, consider
 $$
 \eta_\sigma(h)=\sup _{t \in R^{+}} \int_t^{t+h}\|\sigma(s)\|_{\mathbb{E}_1}^p d s<+\infty.
 $$
 
 Similar to Theorem \ref{th2.4}, we have the following existence and uniqueness result.
 \begin{Theorem}\label{th3.10}
 	Define the set $\Sigma=\left[G_0(t+h) \mid h \in R^{+}\right]_{\mathbb{E}}$ where $G_0 \in \mathbb{E}$ is an arbitrary but fixed symbol function. For any $G \in \Sigma$ and any initial condition $\left(u_\tau, v_{\tau}, \theta_\tau\right) \in H_c, \tau \geq 0$, there exists a unique global solution $\left(u(t), u_t(t), \theta(t)\right) \in \mathcal{H}$, which generates a unique process $\left\{U_G(t, \tau)\right\}(t \geq \tau, \tau \geq 0)$ on $H_c$ of a two-parameter family of operators, such that for any $t \geq \tau, \tau \geq 0$,
 	$$
 	\begin{aligned}
 		& U_G(t, \tau)\left(u_\tau, v_{ \tau}, \theta_\tau\right)=\left(u(t), u_t(t),\theta(t)\right) \in H_c, \\
 		& u(t) \in C\left(R_\tau, W_c\right), \\
 		& u_t(t) \in C\left(R_\tau, E_c\right), \\
 		& \theta(t) \in C\left(R_\tau, L^2(\Omega)\right) .
 	\end{aligned}
 	$$
 \end{Theorem}
  
\begin{Theorem}\label{th3.11}
	 Under the assumption of (\ref{eq3.35}), the family of processes $\left\{U_G(t, \tau)\right\}(G \in \Sigma, t \geq \tau, \tau \geq 0)$, corresponding to (\ref{eq3.33})-(\ref{eq3.34}) has a bounded uniformly (w.r.t. $G \in \Sigma$ ) absorbing set $B_0$ in $H_c$.
	 
\rm\textbf{Proof.} 
 As established in \cite{VV3}, $\mathcal{L}(t)$ is equivalent to $E(t)$, that is, $\mathcal{L}(t)\sim E(t)$. Then according to Lemma \ref{le3.9}, we derive that
 \begin{equation}\label{eq3.36}
 	\frac{d}{dt}E(t)\leq -\xi_1 E(t)+\tilde{C}_1\|g\|^2+\tilde{M}_1.
 \end{equation}
Using assumption (\ref{eq2.6}) and the Gronwall's inequality, we know
\begin{equation}\label{eq3.37}
	\begin{aligned}
		E(t)&\leq E(\tau)e^{-\xi_1(t-\tau)}+(1+\xi_1^{-1})(\tilde{M}_1+\|g\|^2)\\
		&\leq E(\tau)e^{-\xi_1(t-\tau)}+(1+\xi_1^{-1})(\tilde{M}_1+\|g\|_{L_b^2}^2).
	\end{aligned}
\end{equation}

 From Theorem \ref{th3.4}, we infer
 $$
 \|g\|_{L_b^2}^2 \leqslant\left\|g_0\right\|_{L_b^2}^2, \quad \text { for all } g \in \Sigma.
 $$
 
 Therefore, we obtain the uniformly (w.r.t. $G \in \Sigma$ ) absorbing set $B_0$ in $H_c$
 $$
 B_0=\left\{\left(u, u_t,\theta\right) \mid\|\left(u, u_t,\theta\right)\|_{H_c}^2 \leqslant \rho_0^2\right\},
 $$
 where $\rho_0=2\left(1+\xi_1^{-1}\right)\left((\tilde{M}_1+
 \|g_0\|_{L_b^2}^2)\right)$. In other words, for any bounded subset $B$ in $H_c$, there exists a time $t_0=$ $t_0(\tau, B) \geqslant \tau$ such that
 $$
 \bigcup_{G \in \Sigma} U_G(t, \tau) B \subset B_0, \quad \forall t \geqslant t_0.
 $$
$\hfill\qedsymbol$ 
\end{Theorem}
\section{{Uniformly (w.r.t. $G \in \Sigma$ ) asymptotic compactness in $H_c$} }

\quad In this section, we begin with essential preliminaries and subsequently derive a priori energy estimates, drawing on concepts from \cite{II7,AK9,CS8}. We conclude by applying Theorem \ref{th4.2} to establish uniform (w.r.t. $G \in \Sigma)$ asymptotic compactness in $H_c$.

From this point onward, we will consistently denote by $B_0$ the bounded uniformly (w.r.t. $G \in \Sigma$ ) absorbing set obtained in Theorem \ref{th3.11}.

\textbf{\large {4.1. Preliminaries}}
\begin{Definition}\label{de4.1}(\cite{IC4,CS8})
 Let $X$ be a Banach space, $B$ a bounded subset of $X$ and $\Sigma$ a symbol (or parameter) space. We call a function $\phi(\cdot, \cdot ; \cdot, \cdot)$, defined on $(X \times X) \times(\Sigma \times \Sigma)$, to be a contractive function on $B \times B$ if for any sequence $\left\{x_n\right\}_{n=1}^{\infty} \subset B$ and any $\left\{\sigma_n\right\} \subset \Sigma$, there is a subsequence $\left\{x_{n_k}\right\}_{k=1}^{\infty} \subset\left\{x_n\right\}_{n=1}^{\infty}$ and $\left\{\sigma_{n_k}\right\}_{k=1}^{\infty} \subset\left\{\sigma_n\right\}_{n=1}^{\infty}$ such that
	$$
	\lim _{k \rightarrow \infty} \lim _{l \rightarrow \infty} \phi\left(x_{n_k}, x_{n_l} ; \sigma_{n_k}, \sigma_{n_l}\right)=0 .
	$$	
\end{Definition}

We denote the set of all contractive functions on $B \times B$ by $\operatorname{Contr}(B, \Sigma)$.

\begin{Theorem}\label{th4.2}(\cite{IC4})
	Let $\left\{U_\sigma(t, \tau)\right\}, \sigma \in \Sigma$, be a family of processes which satisfies the translation identity in Definition \ref{de3.1} on Banach space X and has a bounded uniformly (w.r.t. $\sigma \in \Sigma$ ) absorbing set $B_0 \subset X$. Moreover, assume that for any $\varepsilon>0$, there exist $T=T\left(B_0, \varepsilon\right)$ and $\phi_T \in \operatorname{Contr}\left(B_0, \Sigma\right)$ such that
	
	$$
	\left\|U_{\sigma_1}(T, 0) x-U_{\sigma_2}(T, 0) y\right\| \leqslant \varepsilon+\phi_T\left(x, y ; \sigma_1, \sigma_2\right), \quad \forall x, y \in B_0, \forall \sigma_1, \sigma_2 \in \Sigma .
	$$
	Then $\left\{U_\sigma(t, \tau)\right\}, \sigma \in \Sigma$, is uniformly (w.r.t. $\sigma \in \Sigma$ ) asymptotically compact in $X$.
\end{Theorem}

\textbf{\large {4.2. A priori estimate}}

Next, we aim to obtain the uniformly $(w.r.t.\,G\in\Sigma)$ asymptotic compactness by establishing the inequality (\ref{eq4.4}). Without loss of generality, we will focus on the strong solutions in the sequence, as the case for weak solutions can be readily addressed using a density argument.

For any $\left(u_\tau^i, v_{ \tau}^i,  \theta_\tau^i\right) \in B_0$, let $\left(u_i(t), v_{i }(t), \theta_i(t)\right)$ be the corresponding solution to $G_i \in \Sigma, i=1,2$.

Let
$$
Z(t)=(z(t)-\varphi(t))^T=\left(u_1(t)-u_2(t), \theta_1(t)-\theta_2(t)\right)^T ,
$$
 then the difference $Z=\left(z, z_t, \varphi\right)$ solves the following problem
\begin{equation}\label{eq4.1}
	\left\{\begin{array}{l}
		z_{tt}-\Delta_e z+\alpha(t) \nabla \varphi+f\left(u^1\right)-f\left(u^2\right)=\nabla \widetilde{p}, \\
		\varphi_t-\kappa(t) \Delta \varphi+\alpha(t) \operatorname{div}  z_t=g_1(x,t)-g_2(x,t),
	\end{array}\right.
\end{equation}
where $\Delta \widetilde{p}=0$,  with initial and boundary conditions
\begin{equation}\label{eq4.2}
	\left\{\begin{array}{l}
		z(\tau)=u^1(\tau)-u^2(\tau), \quad z_t(\tau)= u^1_t(\tau)- u^2_t(\tau), \\
		z(x, t)=0, \varphi(x, t)=0, x \in \partial \Omega, t \geq \tau.
	\end{array}\right.
\end{equation}

We define the functional
\begin{equation}\label{eq4.3}
	E_Z(t):=\frac{1}{2} \int_{\Omega}\left[\left| z_t(t)\right|^2+(2 \mu+\lambda)|\nabla z(t)|^2+|\varphi(t)|^2\right] \mathrm{d} x .
\end{equation}

\begin{Lemma}\label{le4.3}
	Assume the condition (\ref{eq2.6}) is satisfied. Then for any fixed $T>0$, there exist a constant $C_T$ and a function $\phi_T=\phi_T\left(\left(u_0^1, v_0^1,\theta_0^1\right),\left(u_0^2, v_0^2,\theta_0^2\right) ; G_1, G_2\right)$ such that
	\begin{equation}\label{eq4.4}
		\left\|u_1(T)-u_2(T)\right\|_{H_c} \leqslant C_T+\phi_T\left(\left(u_0^1, v_0^1,\theta_0^1\right),\left(u_0^2, v_0^2,\theta_0^2\right) ; G_1, G_2\right),
	\end{equation}
	where $G_i=(0,0,g_i),i=1,2$, $C_T$ and $\phi_T$ depends on $T$.
	
 \rm\textbf{Proof.} Multiplying the equation $(\ref{eq4.1})_1$ by $z_t$ and $(\ref{eq4.1})_2$ by $\varphi$, then summing these products, and utilize the orthogonality of $z_t$ and $\nabla \widetilde{p}$, along with assumption (\ref{eq2.2}), we deduce that 
 \begin{equation}\label{eq4.5}
 	\begin{aligned}
 		\frac{d}{dt}E_Z(t)&=-\kappa(t)\int_\Omega|\nabla\varphi|^2\mathrm{d} x+\int_\Omega(f(u^2)-f(u^1))\cdot z_t\mathrm{d} x+\int_\Omega(g_1-g_2)\varphi \mathrm{d} x\\
 		&\leq -k_0(t)\int_\Omega|\nabla\varphi|^2\mathrm{d} x+\int_\Omega(f(u^2)-f(u^1))\cdot z_t\mathrm{d} x+\int_\Omega(g_1-g_2)\varphi \mathrm{d} x.
 	\end{aligned}
 \end{equation}
 
Intergrating (\ref{eq4.5}) over [ $\sigma, T]$, we get
\begin{equation}\label{eq4.6}
	\begin{aligned}
		E_Z(T)&\leq -k_0\int_\sigma^T\int_\Omega|\nabla\varphi|^2\mathrm{d} xds+\int_\sigma^T \int_{\Omega}\left(f(u^1)-f(u^2)\right)\cdot z_td x d s \\
		&\quad+\int_\sigma^T \int_{\Omega}\left(g_1-g_2\right) \varphi d x d s+E_Z(\sigma).
	\end{aligned}		
\end{equation}

Then intergrating (\ref{eq4.6}) over $[0, T]$, we  obtain 
\begin{equation}\label{eq4.7}
	\begin{aligned}
		E_Z(T)&\leq -k_0\int_{0}^{T}\int_\sigma^T\int_\Omega|\nabla\varphi|^2\mathrm{d} xdsd\sigma+\int_{0}^{T}\int_\sigma^T \int_{\Omega}\left(f(u^1)-f(u^2)\right)\cdot z_td x d sd\sigma \\
		&\quad+\int_{0}^{T}\int_\sigma^T \int_{\Omega}\left(g_1-g_2\right) \varphi d x d sd\sigma+\int_{0}^{T}E_Z(\sigma)d\sigma.
	\end{aligned}
\end{equation}

 By Young's inequality, we derive
 \begin{equation}\label{eq4.8}
  \int_{0}^{T}\int_\sigma^T \int_{\Omega}\left(g_1-g_2\right) \varphi d x d sd\sigma\leq \frac{1}{2} \int_0^T \int_\sigma^T \int_{\Omega}\left|g_1-g_2\right|^2 d x d s d \sigma+\frac{1}{2} \int_0^T \int_\sigma^T \int_{\Omega} |\varphi|^2 d x d s d \sigma.
 \end{equation}

Based on Lemma \ref{le2.1}, the Young's inequality, H$\ddot{o}$lder inequality\cite{YQ30} with $\frac{\rho-1}{2 \rho}+\frac{1}{2 \rho}+\frac{1}{2}=1$, and embedding $\left[H_0^1(\Omega)\right]^N \hookrightarrow\left[L^{2 \rho}(\Omega)\right]^N$, we deduce
\begin{equation}\label{eq4.9}
	\begin{aligned}
		\int_0^T& \int_\sigma^T\int_{\Omega}\left(f\left(u^2\right)-f\left(u^1\right)\right) \cdot  z_t d xd s d \sigma\\  &=\int_0^T \int_\sigma^T\sum_{i=1}^N \int_{\Omega}\left(f_i\left(u^2\right)-f_i\left(u^1\right)\right) \partial_t z_i \mathrm{~d} xd s d \sigma \\
		& \leq 2^{\rho-1} N\int_0^T \int_\sigma^T \sum_{i=1}^N \int_{\Omega}\left(1+\sum_{i=1}^N\left|u_i^2\right|^{\rho-1}+\sum_{i=1}^N\left|u_i^1\right|^{\rho-1}\right)\left|z_i \| \partial_t z_i\right| d xd s d \sigma \\
		& \leq C\int_0^T \int_\sigma^T\left(1+\left\|u^1\right\|_{2 \rho}^{\rho-1}+\left\|u^2\right\|_{2 \rho}^{\rho-1}\right)\|z\|_{2 \rho}\left\| z_t\right\|d s d \sigma \\
		& \leq C\int_0^T \int_\sigma^T\left(1+\left\|u^1\right\|_{\left[H_0^1(\Omega)\right]^N}^{\rho-1}+\left\|u^2\right\|_{\left[H_0^1(\Omega)\right]^N}^{\rho-1}\right)\|z\|_{2 \rho}\left\| z_t\right\|d s d \sigma \\
		& \leq C_B\int_0^T \int_\sigma^T\|z\|_{2 \rho}\left\| z_t\right\|d s d \sigma \\
		& \leq \frac{C_B^2}{2}\int_0^T \int_\sigma^T\|z\|_{2 \rho}^2d s d \sigma+\frac{1}{2} \int_0^T \int_\sigma^T\int_{\Omega}\left| z_t\right|^2 \mathrm{~d} xd s d \sigma.
	\end{aligned}
\end{equation}

According to (\ref{eq4.3}), we conclude
\begin{equation}\label{eq4.10}
	\frac{1}{2} \int_0^T \int_\sigma^T \int_{\Omega} |\varphi|^2 d x d s d \sigma+\frac{1}{2} \int_0^T \int_\sigma^T\int_{\Omega}\left| z_t\right|^2 \mathrm{~d} xd s d \sigma\leq \int_0^T \int_\sigma^TE_Z(\sigma)d s d \sigma.
\end{equation}

Using (\ref{eq3.37}), we have
\begin{equation}\label{eq4.11}
	\int_{0}^{T}E_Z(\sigma)d\sigma\leq \int_{0}^{T}\left[E(\tau)e^{-\xi_1(t-\tau)}+(1+\xi_1^{-1})\left(\tilde{M}_1+\|g_0\|_{L_b^2}^2\right)\right]d\sigma.
\end{equation}

Hence, there exists a positive constant 
 $C_M=C(T, \tau, \mu)>0$ such that
\begin{equation}\label{eq4.12}
	\int_0^T \int_\sigma^T E_Z(\sigma) d s d \sigma+\int_0^T E_Z(\sigma) d \sigma \leq C_M .
\end{equation}

Inserting (\ref{eq4.8})-(\ref{eq4.12}) in 
(\ref{eq4.7}), we conclude
\begin{equation}\label{eq4.13}
	T E_Z(T) \leq C_M+\frac{C_B^2}{2}\int_0^T \int_\sigma^T\|u^1-u^2\|_{2 \rho}^2d s d \sigma+\frac{1}{2} \int_0^T \int_\sigma^T \int_{\Omega}\left|g_1-g_2\right|^2 d x d s d \sigma .
\end{equation}

Set
\begin{equation}\label{eq4.14}
	\begin{aligned}
		\phi\left(\left(u_0^1, v_{0}^1,  \theta_0^1\right),\left(u_0^2, v_{0}^2,  \theta_0^2\right) ; G_1, G_2\right)&=\frac{C_B^2}{2}\int_0^T \int_\sigma^T\|u^1-u^2\|_{2 \rho}^2d s d \sigma\\
		&\quad+\frac{1}{2} \int_0^T \int_\sigma^T \int_{\Omega}\left|g_1-g_2\right|^2 d x d s d \sigma.
	\end{aligned}	
\end{equation}

Consequently, we get
\begin{equation}\label{eq4.15}
	E_Z(T) \leqslant \frac{C_M}{T}+\frac{1}{T} \phi_T\left(\left(u_0^1, v_0^1,\theta_0^1\right),\left(u_0^2, v_0^2,\theta_0^2\right) ; G_1, G_2\right) .
\end{equation}

The proof is complete.
$\hfill\qedsymbol$ 
\end{Lemma}

\textbf{\large {4.3. Uniformly asymptotic compactness}}

In this subsection, we shall prove the uniformly (w.r.t. $G \in \Sigma$ ) asymptotic compactness in $H_c$, which is given in the following theorem.
\begin{Theorem}\label{th4.4}
	 Assume that $G$ satisfies (\ref{eq3.35}) and $g$ satisfies (\ref{eq2.6}), then the family of processes $\left\{U_G(t, \tau)\right\}, G \in \Sigma$, corresponding to (\ref{eq3.33})-(\ref{eq3.34}), is uniformly (w.r.t. $G \in \Sigma$ ) asymptotically compact in $H_c$.
	
 \rm\textbf{Proof.} Given that the family of processes $\left\{U_G(t, \tau)\right\}, G \in \Sigma$, has a bounded uniformly absorbing set. According to Lemma \ref{le4.3}, for any fixed $\varepsilon>0$, we can select a sufficiently large $T$, so that
	$$
	\frac{C_M}{T} \leqslant \varepsilon.
	$$
	
	Consequently, by Theorem \ref{th4.2}, it is sufficient to prove $\phi_T(\cdot, \cdot ; \cdot, \cdot) \in \operatorname{Contr}\left(B_0, \Sigma\right)$ for each fixed $T$.
	From the proof procedure of Theorem \ref{th3.11}, we infer that for any fixed $T$,
	\begin{equation}\label{eq4.16}
		\bigcup_{G \in \Sigma} \bigcup_{t \in[0, T]} U_G(t, 0) \quad B_0 \text { is bounded in } E_0,
	\end{equation}
	and the bound depends on $T$.
	
	Consider the sequence of solutions $\left(u_n, u_{n_t},\theta_n\right)$ corresponding to the initial data $\left(u_0^n, v_0^n,\theta_0^n\right) \in B_0$ associated with the symbols $G_n \in \Sigma$, $n=1,2, \ldots$. Based on equation  (\ref{eq4.16}), and without loss of generality (possibly by considering a subsequence), we assume that
	\begin{equation}\label{eq4.17}
			\begin{aligned}
			& u_n \rightarrow u \quad \star \text {-weakly in } L^{\infty}\left(0, T ; (H_0^1(\Omega))^N\right), \\
			& u_{n_t} \rightarrow u_t \quad \star \text {-weakly in } L^{\infty}\left(0, T ; (L^2(\Omega))^N\right), \\
			& u_n \rightarrow u \quad \text { in } L^2\left(0, T ; (L^2(\Omega))^N\right), \\
			& u_n \rightarrow u \quad \text { in } L^k\left(0, T ; (L^k(\Omega))^N\right), \\
			& u_n(0) \rightarrow u(0) \quad \text { and } \quad u_n(T) \rightarrow u(T) \quad \text { in } L^k(\Omega),
		\end{aligned}
	\end{equation}
	for $k<\frac{N}{N-2}$, where we have used the compact embedding $\left[H_0^1(\Omega)\right]^N \hookrightarrow\left[L^{2 \rho}(\Omega)\right]^N$.
	
	From (\ref {eq2.6}), we can infer
	\begin{equation}\label{eq4.18}
		\lim _{n \rightarrow+\infty} \lim _{m \rightarrow+\infty} \int_0^T \int_\sigma^T \int_{\Omega}\left|g_n(x, s)-g_m(x, s)\right|^2 d x d s d \sigma=0.
	\end{equation}
	
	From (\ref{eq4.17}) and the compact embedding $\left[H_0^1(\Omega)\right]^N \hookrightarrow\left[L^{2 \rho}(\Omega)\right]^N$, by the Aubin-Lions theorem, there exists a subsequence $\left\{u^{n}\right\}$ such that converges for some $u$ in $L^2\left(0,T;\left[L^{2 \rho}(\Omega)\right]^N\right)$. Thus, we have
	\begin{equation}\label{eq4.19}
		\lim _{n \rightarrow+\infty} \lim _{m \rightarrow+\infty} \int_0^T \int_\sigma^T\|u^n-u^m\|_{2 \rho}^2d s d \sigma=0.
	\end{equation}
	
Finally, combining (\ref{eq4.18})-(\ref{eq4.19}), we get $\phi_T(\cdot, \cdot ; \cdot, \cdot) \in \operatorname{Contr}\left(B_0, \Sigma\right)$  immediately.	
$\hfill\qedsymbol$ 	
\end{Theorem}

\textbf{\large {4.4. Existence of uniform attractor}}
\begin{Theorem}\label{th4.5}
Under the assumption (\ref{eq2.1})-(\ref{eq2.6}) and the definition of $\Sigma$, then the family of processes $\left\{U_G(t, \tau)\right\}(G \in \Sigma, t \geq$ $\tau, \tau \geq 0$ ) corresponding to (\ref{eq2.7})-(\ref{eq2.8}) has a compact uniform (w.r.t. $G \in \Sigma$ ) attractor $\mathcal{A}_{\Sigma}$.
	
\rm\textbf{Proof.} Theorems \ref{th3.11} and \ref{th4.4} imply the existence of a uniform attractor immediately.
$\hfill\qedsymbol$ 	
\end{Theorem}

$\mathbf{Acknowledgment}$

This work was supported by the TianYuan Special Funds of the NNSF of China with contract number 12226403, the NNSF of China with contract
No.12171082, the fundamental research funds for the central universities with contract numbers 2232022G-13, 2232023G-13, 2232024G-13.

\end{document}